\documentclass[11pt,leqno]{article}

\usepackage{amsmath,amsfonts,amscd,amssymb,theorem}

\long\def\comment#1\endcomment{}

\makeatletter
\begingroup
\gdef\th@dotted{\normalfont\itshape
  \def\@begintheorem##1##2{%
        \item[\hskip\labelsep \theorem@headerfont ##1\ ##2.]}%
\def\@opargbegintheorem##1##2##3{%
   \item[\hskip\labelsep \theorem@headerfont ##1\ ##2\ (##3).]}}
\endgroup
\makeatother

\theoremstyle{dotted}

\newtheorem{theorem}{Theorem}[section]
\newtheorem{lemma}[theorem]{Lemma}
\newtheorem{conj}[theorem]{Conjecture}
\newtheorem{prop}[theorem]{Proposition}
\newtheorem{corr}[theorem]{Corollary}

\makeatletter
\begingroup
\gdef\th@upshape{\normalfont
  \def\@begintheorem##1##2{%
        \item[\hskip\labelsep \theorem@headerfont ##1\ ##2.]}%
\def\@opargbegintheorem##1##2##3{%
   \item[\hskip\labelsep \theorem@headerfont ##1\ ##2\ (##3).]}}
\endgroup
\makeatother

\theoremstyle{upshape}

\newtheorem{defn}[theorem]{Definition}
\newtheorem{remark}[theorem]{Remark}

\makeatletter
\renewcommand{\subsection}{\@startsection{subsection}{2}{0pt}{-3ex
plus -1ex minus -0.2ex}{-2mm plus -0pt minus
-2pt}{\normalfont\bfseries}} \makeatother

\makeatletter
\@addtoreset{equation}{section}
\makeatother

\newcommand{\cntrct}                
{\hspace{2pt}\raisebox{1pt}{\text{$\lrcorner$}}\hspace{2pt}}

\newcommand{\proof}[1][Proof.]{\smallskip\noindent{\em #1}}
\def\endproof{\hfill\ensuremath{\square}\par\medskip}

\def\eqref#1{\thetag{\ref{#1}}}

\let\latexref=\ref
\def\ref#1{{\normalfont{\latexref{#1}}}}

\newcommand{\wt}{\widetilde}
\newcommand{\wh}{\widehat}

\newcommand{\idot}{{\:\raisebox{1pt}{\text{\circle*{1.5}}}}}
\newcommand{\hdot}{{\:\raisebox{3pt}{\text{\circle*{1.5}}}}}
\newcommand{\calo}{{\cal O}}

\newcommand{\thebull}{\thetag{$\bullet$}}

\newcommand{\gm}{{\mathbb{G}_m}}

\newcommand{\6}{\partial}

\newcommand{\C}{{\mathbb C}}
\newcommand{\Z}{{\mathbb Z}}
\newcommand{\R}{{\mathbb R}}

\newcommand{\X}{{\mathfrak X}}
\newcommand{\Y}{{\cal Y}}
\newcommand{\T}{{\cal T}}
\newcommand{\N}{{\cal N}}
\newcommand{\m}{{\mathfrak m}}
\newcommand{\g}{{\mathfrak g}}

\newcommand{\id}{\operatorname{\sf id}}
\renewcommand{\dim}{\operatorname{\sf dim}}
\newcommand{\codim}{\operatorname{\sf codim}}
\newcommand{\rk}{\operatorname{\sf rk}}
\newcommand{\cchar}{\operatorname{\sf char}}
\newcommand{\Spec}{\operatorname{Spec}}
\newcommand{\Ker}{\operatorname{{\sf Ker}}}
\newcommand{\End}{\operatorname{End}}
\newcommand{\red}{\operatorname{\sf red}}
\newcommand{\ad}{\operatorname{\sf ad}}
\newcommand{\gr}{\operatorname{\sf gr}}

\newcommand{\whotimes}{\operatorname{\wh{\otimes}}}

\title{Symplectic singularities from\\ the Poisson point of view}

\author{D. Kaledin\thanks{Partially supported by CRDF grant
RM1-2354-MO02.}}

\begin{document}

\maketitle

\tableofcontents

\section*{Introduction}

In symplectic geometry, it is often useful to consider the so-called
{\em Poisson bracket} on the algebra of functions on a $C^\infty$
symplectic manifold $M$. The bracket determines, and is determined
by, the symplectic form; however, many of the features of symplectic
geometry are more conveniently described in terms of the Poisson
bracket. When one turns to the study of symplectic manifolds in the
holomorphic or algebro-geometric setting, one expects the Poisson
bracket to be even more useful because of the following observation:
the bracket is a purely algebraic structure, and it generalizes
immediately to {\em singular} algebraic varieties and
complex-analytic spaces.

The appropriate notion of singularities for symplectic algebraic
varieties has been introduced recently by A. Beauville \cite{B} and
studied by Y. Namikawa \cite{N1}, \cite{N2}. However, the theory of
singular symplectic algebraic varieties is in its starting stages;
in particular, to the best of our knowledge, the Poisson methods
have not been used yet. This is the goal of the present paper.

\bigskip

Our results are twofold. Firstly, we prove a simple but useful
structure theorem about symplectic varieties (Theorem~\ref{th.1})
which says, roughly, that any symplectic variety admits a canonical
stratification with a finite number of symplectic strata (in the
Poisson language, a symplectic variety considered as a Poisson space
has a finite number of symplectic leaves). In addition, we show
that, locally near a stratum, the variety in question admits a nice
decomposition into the product of the stratum itself and a
transversal slice. Secondly, we study natural group actions on a
symplectic variety and we prove that, again locally, a symplectic
variety always admits a non-trivial action of the one-dimensional
torus $\gm$ (Theorem~\ref{th.2}). This is a rather strong
restriction on the type of singularities a symplectic variety might
have.

Unfortunately, the paper is much more eclectic than we would like.
Moreover, one of the two main results is seriously flawed: we were
not able to show that the $\gm$-action provided by
Theorem~\ref{th.2} has positive weights. However, all the results
has been known to the author for a couple of years now, and it seems
that any improvement would require substantially new methods. Thus
we have decided to publish the statements ``as is''.

In addition, we separately consider a special (and relatively rare)
situation when a symplectic variety admits a crepant resolution of
singularities. We prove that the geometry of such a resolution is
very restricted: it is always semismall, and the Hodge structure on
the cohomology of its fibers is pure and Hodge-Tate.

\bigskip

Our approach, for better or for worse, is to try to use Poisson
algebraic methods as much as possible, getting rid of actual
geometry at an early stage. The paper is organized as follows. In
the first section we recall all the necessary definitions, both from
the Poisson side of the story and from the theory of symplectic
singularities. We also introduce two particular classes of Poisson
schemes which we call {\em holonomic} and {\em locally exact}. In
the second section we show that symplectic varieties give examples
of Poisson schemes lying in both of these classes. The main
technical tool here is the beautiful canonical resolution of
singularities discovered in the last two decades (see, for example,
\cite{BM}). Then in Subsection~\ref{crep}, we study the geometry of
crepant resolutions. The remainder of the paper is purely
algebraic. First, we prove the structure theorems for holonomic
Poisson schemes. Then we show that is a scheme is in addition
locally exact, then it admits, again locally, a $\gm$-action.

\subsection*{Acknowledgements.} In the recent years, I have benefitted
from many opportunities to discuss symplectic singularities with
colleages; the discussion were very fruitful and often led to better
understanding on my part. I would like to thank, in particular,
R. Bezrukavnikov, H. Esnault, P. Etingof, B. Fu, V. Ginzburg,
D. Huybrechts, D. Kazhdan, A. Kuznetsov, M. Lehn, Y. Namikawa,
M. Verbitsky, E. Viehweg, J. Wierzba -- and I would like to
apologize to those I forgot to mention here. The results were
presented at an Oberwolfach meeting on singularities in 2003. I am
grateful to the organizers, in particular to D. van Straten, for
inviting me and giving me an opportunity to benefit from the
beautiful atmosphere of this great institution. I am also grateful
to A. Bondal, B. Fu, V. Ginzburg and Y. Namikawa who have read the
first version of this paper, suggested several improvements and
found some gaps in the proofs. Finally, I would like to thank the
referee for a thoughtful report and important suggestions and
corrections.

\section{Generalities on Poisson schemes.}

Fix once and for all a base field $k$ of characteristic $\cchar k =
0$.

\begin{defn}\label{poi.def}
A {\em Poisson algebra} over the field $k$ is a commutative algebra
$A$ over $k$ equipped with an additional skew-linear operation
$\{-,-\}:A \otimes A \to A$ such that
\begin{equation}\label{poi}
\begin{aligned}
\{a,bc\} &= \{a,b\}c + \{a,c\}b, \\
0 &= \{a,\{b,c\}\} + \{b,\{c,a\}\} + \{c,\{a,b\}\},
\end{aligned}
\end{equation}
for all $a,b,c \in A$. An ideal $I \subset A$ is called a {\em
Poisson ideal} if $\{i,a\} \in I$ for any $i \in I$, $a \in A$.
\end{defn}

We will always assume that a Poisson algebra $A$ has a unit element
$1 \in A$ such that $\{1,a\}=0$ for every $a \in A$.

\begin{defn}
A {\em Poisson scheme} over $k$ is a scheme $X$ over $k$ equipped
with a skew-linear bracket in the structure sheaf $\calo_X$
satisfying \eqref{poi}.
\end{defn}

If $A$ is a Poisson algebra over $k$, then $X = \Spec A$ is a
Poisson scheme. The reduction, every irreducible component, any
completion and the normalization of a Poisson scheme are again
Poisson schemes (\cite{K1}). We will say that a Poisson scheme is
{\em local} if it is the spectrum of a local Poisson algebra $A$
whose maximal ideal $\m$ is a Poisson ideal.

Let $X$ be a Poisson scheme. For every local function $f$ on $X$,
the bracket $\{f,-\}$ is by definition a derivation of the algebra
of functions, hence a vector field of $X$, denoted by $H_f$. Vector
fields of the form $H_f$ are called {\em Hamiltonian} vector
fields. Moreover, the Poisson bracket is a derivation with respect
to each of the two arguments. Therefore it can be expressed as
\begin{equation}\label{theta}
\{f,g\} = \Theta(df \wedge dg),
\end{equation}
where
$$
\Theta:\Lambda^2\Omega^1_X \to \calo_X
$$
is an $\calo_X$-linear map. The map $\Theta$ is called the {\em
Poisson bivector}. By the Jacobi identity part of \eqref{poi}, we
have $H_f(\Theta)=0$ for every Hamiltonian vector field $H_f$
(Hamiltonian vector fields preserve the Poisson bivector). If $X$ is
smooth -- for instance, if it is a point -- then the cotangent sheaf
$\Omega^1_X$ is flat and $\Theta$ gives a skew-symmetric bilinear
form on this sheaf.

Given a closed subscheme $Y \subset X$, we will say that $Y$ is a
{\em Poisson subscheme} if it is locally defined by a Poisson ideal
in $\calo_X$. Equivalently, a subscheme is Poisson if locally it is
preserved by all Hamiltonian vector fields (in other words, all
Hamiltonian vector fields are tangent to $Y$). In this case, $Y$
inherits the structure of a Poisson scheme.

\begin{defn}
\begin{enumerate}
\item An Noetherian intergal Poisson scheme $X$ over $k$ with
generic point $\eta \in X$ is called {\em generically
non-degenerate} if the Poisson bivector $\Theta$ gives a
non-degenerated skew-symmetric form on the cotangent module
$\Omega^1(\eta/k)$.
\item A Noetherian integral Poisson scheme $X$ is called {\em
holonomic} if every integral Poisson subscheme $Y \subset X$ is
generically non-degenerate.
\end{enumerate}
\end{defn}

In particular, a holonomic Poisson scheme $X$ is itself generically
non-degenerate. Moreover, every integral Poisson subscheme $Y
\subset X$ of a holonomic Poisson scheme $X$ is obviously
holonomic. By definition, $X$ itself and every such subscheme $Y
\subset X$ must be even-dimensional over $k$. The normalization $X'$
of a holonomic Poisson scheme is also holonomic (for every prime
ideal $J' \subset A'$ in the normalization $A'$ of a Poisson algebra
$A$, the intersection $J = J' \cap A \subset A$ is a prime ideal in
$A$, $A'/J'$ is generically \'etale over $A/J$, and if $J'$ is
Poisson, $J$ is obviously also Poisson).

The notion of a holonomic Poisson scheme has any meaning only for
singular Poisson schemes; for a smooth scheme $X$ it is vacuous
because of the following.

\begin{lemma}\label{sm}
Let $X$ be a smooth Poisson scheme over $k$. Assume that $X$ is
holonomic. Then the Poisson bivector $\Theta$ is non-degenerate
everywhere on $X$, and the only Poisson subscheme in $X$ is $X$
itself.
\end{lemma}

\proof{} Let $2n=\dim X$. The top degree power $\Theta^n$ of the
Poisson bivector $\Theta$ is a section of the anticanonical bundle
$K^{-1}_X$. Moreover, $\Theta$ is non-degenerate if and only if
$\Theta^n$ is non-zero. Let $D \subset X$ be the zero locus of
$\Theta^n$. It is either empty, or a divisor in $X$. All the
Hamiltonian vector fields preserve $\Theta$, hence also $\Theta^n$
and $D$. Thus $D \subset X$ is a Poisson subscheme. But since $X$ is
by assumption holonomic, $D$ must be even-dimensional -- in
particular, it cannot be a divisor. This proves the first claim. To
prove the second, let $Y \subset X$ be a Poisson subscheme, and let
$y \in Y$ be a closed point in the smooth part of $Y$. Then all
Hamiltonian vector fields $H_f$ are by definition tangent to $Y$ at
$y$. But since $\Theta$ is non-degenerate, Hamiltonian vector fields
span the whole tangent space $T_yX$, and we have $Y = X$.  \endproof

Holonomic Poisson schemes are the first special class of Poisson
schemes that we will need in this paper. To introduce the second
class, we give the following definition.

\begin{defn}
A Poisson algebra $A$ over $k$ is called {\em exact} if there exists
a derivation $\xi:A \to A$ such that
\begin{equation}\label{dilat}
\xi(\{a,b\}) = \{a,b\} + \{\xi(a),b\} + \{a,\xi(b)\}
\end{equation}
for any $a,b \in A$.
\end{defn}

\begin{remark}
This definition is motivated by the theory of {\em Poisson
cohomology}, see e.g. \cite{KG}. The Poisson bivector $\Theta$
defines a degree-$2$ Poisson cocycle on $A$, while any derivation
$\xi:A \to A$ defines a degree-$1$ Poisson cochain. Equation
\eqref{dilat} then says that $\Theta = \delta(\xi)$, where $\delta$
is the Poisson differential. We will not need this, so we do not
give any details and refer the interested reader to
\cite[Appendix]{KG}.
\end{remark}

\begin{defn}
A Poisson scheme $X$ of finite type over $k$ is called {\em locally
exact} if for any closed point $x \in X$, the completed local ring
$\wh{\calo}_{X,x}$ is an exact Poisson $k$-algebra.
\end{defn}

We will see in Section~\ref{str} that for holonomic Poisson schemes,
local exactness passes to Poisson subschemes.

\section{Symplectic singularities.}

\subsection{Statements.}
In this Section, assume that the base field $k$ is a subfield $k
\subset \C$ of the field of complex numbers. Let $X$ be an algebraic
variety -- that is, an integral scheme of finite type over
$k$. Assume that $X$ is normal and even-di\-men\-si\-on\-al, of
dimension $\dim X = 2n$. Assume given a non-degenerate closed
$2$-form $\Omega \in \Omega^2(U)$ on the smooth open part $U \subset
X$.

\begin{defn}[\cite{B},\cite{N2}]
One says that $X$ is a {\em symplectic variety} -- or, equivalently,
that $X$ has {\em symplectic singularities} -- if the $2$-form
$\Omega$ extends to a (possibly degenerate) $2$-form on a resolution
of singularities $\wt{X} \to X$.
\end{defn}

Note that since for any two smooth birational varieties $X_1$, $X_2$
and any integer $k$ we have 
$$
H^0(X_1,\Omega^k(X_1)) = H^0(X_2,\Omega^k(X_2)),
$$ 
this definition does not depend on the choice of the resolution
$\wt{X}$ (and indeed, one could have said ``any resolution of
singularities'' right in the definition). Symplectic singularities
are always canonical, hence rational (see \cite{B}).

Any normal symplectic variety $X$ is a Poisson scheme. Indeed, since
$X$ is normal, it suffices to define the bracket $\{f,g\}$ of any
two local function $f,g \in \calo_X$ outside of the singular
locus. Thus we may assume that $X$ is smooth. Then, since $\Omega$
is by assumption non-degenerate, it gives an identification $\T(X)
\cong \Omega^1(X)$ between the tangent and the cotangent bundles on
$X$, and this identification in turn gives a bivector $\Theta \in
\Lambda^2\T(X)$. It is well-known that $\Theta$ is a Poisson
bivector for some Poisson structure if and only if $\Omega$ is
closed. We will say that a smooth Poisson scheme $X$ is {\em
symplectic} if the Poisson structure on $X$ is obtained by this
construction from a non-degenerate closed $2$-form $\Omega$. A
smooth Poisson scheme $X$ is symplectic if and only if the Poisson
bivector $\Theta$ is non-degenerate; the symplectic form $\Omega$ is
uniquely defined by $\Theta$. By Lemma~\ref{sm}, this is also
equivalent to saying that the smooth Poisson scheme $X$ is
holonomic.

Exactness for symplectic varieties means exactly what one would
expect.

\begin{lemma}\label{exl}
Let $X = \Spec A$ be a normal affine symplectic variety. Then the
Poisson algebra $A$ is exact if and only if the symplectic form
$\Omega$ satisfies
$$
\Omega = d \alpha
$$
for some $1$-form $\alpha \in \Omega^1(U)$ on the non-singular part
$U \subset X$.
\end{lemma}

\proof{} Using \eqref{theta}, it is easy to check that \eqref{dilat}
in the symplectic case means exactly that
$$
\xi(\Omega) = \Omega,
$$
where $\xi$ acts by Lie derivative. If $A$ is exact, then by the
Cartan homotopy formula we have
$$
\Omega = \xi(\Omega) = \xi \cntrct d\Omega + d(\xi \cntrct \Omega) =
d(\xi \cntrct \Omega).
$$
Conversely, if $\Omega = d\alpha$, then
$$
\xi = \alpha \cntrct \Theta \in \T(U)
$$
satisfies \eqref{dilat}. Since $X$ is normal, $\xi$ extends to a
derivation of the whole algebra $A$.
\endproof

We can now state the main two results of the paper.

\begin{theorem}\label{th.1}
Let $X$ be a normal symplectic variety. Then there exists a finite
stratification $X_i \subset X$ by irreducible Poisson subschemes
such that all the open strata $X_i^o \subset X_i$ are smooth and
sympletic. The only integral Poisson subschemes in $X$ are closed
strata. Moreover, for any closed point $x \in X_i^o \subset X$, the
formal completion $\wh{X}_x$ of $X$ at $x \in X$ admits a product
decomposition
\begin{equation}\label{prod}
\wh{X}_x = \Y_x \times \wh{X^o_i}_x,
\end{equation}
where $\wh{X^o_i}_x$ is the formal completion of the stratum $X_i^o$
at $x$ and $\Y_x$ is a local formal Poisson scheme and a symplectic
variety. The product decomposition is compatible with the Poisson
structures and the symplectic forms.\endproof
\end{theorem}

\begin{theorem}\label{th.2}
Let $X$ be a normal symplectic variety. Then each of the transversal
slices $Y_x$ provided by Theorem~\ref{th.1} admits an action of the
group $\gm$ such that for any point $\lambda \in \gm(k)=k^*$
$$
\lambda \cdot \Omega_Y = \lambda^l\Omega_Y,
$$
where $\Omega_Y$ is the symplectic form on the smooth open part of
$Y_x$, and $l \neq 0$ is some integer.\endproof
\end{theorem}

The product decomposition \eqref{prod} should be understood in the
formal scheme sense (the spectrum of a {\em completed} tensor
product). It is unfortunate that we have to pass to formal
completions in the second part of Theorem~\ref{th.1} and in
Theorem~\ref{th.2}; however, this seems to be inevitable. We do not
know to what extent the product decomposition \eqref{prod} can be
globalized. Note that the action in Theorem~\ref{th.2} must be
non-trivial (so that its existence is a non-trivial restriction on
the geometry of the transversal slice). The strata in
Theorem~\ref{th.1} are what is known as {\em symplectic leaves} of
the Poisson scheme $X$; in particular, we prove that there is only a
finite number of those. Theorem~\ref{th.1} is proved in
Section~\ref{str}, and Theorem~\ref{th.2} is proved in
Section~\ref{grp}. Both are actually direct corollaries of the
corresponding Poisson statements and the following theorem, which is
the main result of this Section.

\begin{theorem}\label{main}
Let $X$ be a normal symplectic variety. Then $X$ is holonomic and
locally exact. Moreover, the normalization of every Poisson
subscheme $Y \subset X$ is also a symplectic variety.
\end{theorem}

Before we turn to the proof of Theorem~\ref{main}, we would like to
note that the converse statement is false, at least in dimension $2$
-- there exist normal Poisson varieties which are holonomic and
locally exact, but not symplectic. In fact, every weakly Gorenstein
normal surface singularity is automatically Poisson and
holonomic. It is often locally exact -- for example, in the case
when it admits a good $\C^*$-action. However, the only symplectic
singularities in dimension $2$ are rational double points. 

On the other hand, if a variety $X$ is non-singular outside of
codimension $4$, then every $2$-form $\Omega$ on the smooth locus $U
\subset X$ extends without poles to any smooth resolution $\wt{X}
\to X$ -- this follows from the beautiful theorem of J. Steenbrink
and D. van Straten \cite{SS}, generalized by H. Flenner
\cite{F}. Thus every holonomic Poisson variety non-singular outside
of codimension $4$ is automatically a symplectic variety.

In practice, holonomic Poisson varieties usually {\em are} singular
in $\codim 2$, but these singularities are canonical -- locally, we
have a product of a smooth scheme and a transversal slice which is
just a Du Val point. Thus in codimension $2$ the singularty is
symplectic. Unfortunately, the general extension theorem of
Flenner-Steenbrink-van Straten says nothing at all about a situation
of this type. There is one partial result, however. It has been
proved by Y. Namikawa \cite[Theorem 4]{N2} that a normal variety $X$
with {\em rational} singularities equipped with a symplectic form on
the smooth part $X^{reg} \subset X$ is automatically a symplectic
variety.

\begin{remark}\label{contact} 
In \cite{CF}, F. Campana and H. Flenner defined the so-called {\em
contact singularities} and proved that isolated contact
singularities do not exist. They also conjectured that every contact
singularity is the product of a symplectic singularity and an affine
line. This conjecture should follow more or less directly from our
Theorem~\ref{th.1} -- one treats a contact singularity of dimension
$2n+1$ as a symplectic singularity of dimension $2n+2$ equipped with
a $\gm$-action, and shows that our product decomposition is
compatible with the $\gm$-action. However, it seems that this
statement is not as interesting as it might have been, because the
notion of contact singularity is too restrictive. Essentially,
Campana and Flenner require that a smooth resolution $\wt{X} \to X$
admits a contact structure with trivial contact line bundle. In our
opinion, especially from the point of view of our
Theorem~\ref{th.2}, it would be more interesting to allow contact
line bundles which are not pulled back from line bundles on $X$.
\end{remark}

\subsection{Geometry of resolutions.}

Let $X$ be a symplectic variety. Recall that $X$ has rational
singularities, so that for every smooth resolution $\pi:\wt{X} \to
X$ we have $R^i\pi_*\calo_{\wt{X}} = 0$ for $i \geq 1$. As usual,
this implies in particular that $R^1\pi_*\Z = 0$ in analytic
topology (consider the exponential exact sequence on $\wt{X}$). The
form $\Omega$ on $\wt{X}$ defines a de Rham cohomology class
$[\Omega] \in H^2_{DR}(\wt{X})$. Recall that if we extend the field
of definition from $k \subset \C$ to $\C$, then by the comparison
theorem, $H^2_{DR}(\wt{X})$ is isomorphic to the topological
cohomology group $H^2(\wt{X},\C)$, computed in analytic topology.

\begin{lemma}\label{from.b}
Let $\pi:\wt{X} \to X$ be a smooth projective resolution of a
symplectic variety $X$ over $\C$, and let $[\Omega] \in
H^2(\wt{X},\C)$ be the cohomology class of the associated $2$-form
$\Omega$ on $\wt{X}$.  Then there exists a cohomology class
$[\Omega_X] \in H^2(X,\C)$ such that $[\Omega] = \pi^*[\Omega_X] \in
H^2(\wt{X},\C)$.
\end{lemma}

\proof{} The proof is an application of a beautiful idea of
J. Wierzba \cite{W}. Consider the Leray spectral seqeunce of the map
$\pi:\wt{X} \to X$ and the associated three-step filtration on
$H^2(\wt{X},\C)$ whose graded pieces are subquotients of
$H^p(X,R^q\pi_*\C)$, $p+q=2$. We have to show that $[\Omega] \in
\pi^*(H^2(X,\C)) = \pi^*(H^2(X,\pi_*\C)) \subset
H^2(\wt{X},\C)$. Since $R^1\pi_*\Z=0$, we also have $R^1\pi_*\C=0$;
therefore the middle term in the associated graded quotient of
$H^2(\wt{X},\C)$ with respect to the Leray filtration vanishes, and
it suffices to prove that the projection $H^2(\wt{X},\C) \to
H^0(X,R^2\pi_*\C)$ annihilates $[\Omega]$. By proper base change,
this is equivalent to proving that for every closed point $x \in X$,
$[\Omega]$ restricts to $0$ on the fiber $Z_x = \pi^{-1}(u) \subset
\wt{X}$ of the map $\pi:\wt{X} \to X$.

Take such a fiber $Z_x$ (with the reduced scheme structure). By
\cite{D}, the cohomology groups $H^\hdot(Z_x,\C)$ carry a natural
mixed Hodge structure. To see it explcitly, one chooses a simplicial
resolution $Z_\idot$ of the scheme $Z_x$ by smooth proper schemes.
Since $\pi:\wt{X} \to X$ is projective, $Z_x$ is also projective,
and we can assume that $Z_\idot$ is a resolution by projective
smooth schemes. Then the Hodge filtration $F^\hdot$ on
$H^\hdot(Z_x,\C) \cong H^\hdot(Z_\idot,\C) = H^\hdot_{DR}(Z_\idot)$
is induced by the stupid filtration on the de Rham complex
$\Omega^\hdot(Z_\idot)$.

Now, denote by $F_X^\hdot$ the filtration induced on
$H^\hdot_{DR}(\wt{X})$ by the stupid filtration on the de Rham
complex $\Omega^\hdot(\wt{X})$. Since $\wt{X}$ is not compact, the
corresponding spectral sequence does not degenerate. Nevertheless,
the restriction map $P:H^\hdot_{DR}(\wt{X}) \to
H^\hdot_{DR}(Z_\idot) \cong H^\hdot(Z_x,\C)$, being induced by the
natural map from $Z_\idot$ to the constant simplicial scheme
$\wt{X}$, obviously sends $F_X^lH^\hdot(\wt{X})$ into
$F^lH^\hdot(Z_x,\C)$.

Since $R^{\geq 1}\pi_*\calo_{\wt{X}}=0$, we may replace $X$ with an
affine neighborhood of the point $x \in X$ and assume that $H^{\geq
1}(\wt{X},\calo_{\wt{X}}) = 0$. Then
$H^l_{DR}(\wt{X})=F^1_XH^l_{DR}(\wt{X})$ for any $l \geq
1$. Applying complex conjugation on $H^l_{DR}(\wt{X}) \cong
H^l(\wt{X},\R) \otimes \C$, we deduce that
$H^l_{DR}(\wt{X})=\overline{F}^1_XH^l_{DR}(\wt{X})$, where
$\overline{F}^\hdot_X$ is the filtration complex-conjugate to
$F^\hdot_X$. In particular, $[\Omega] \in
\overline{F}^1_XH^2_{DR}(\wt{X})$. On the other hand, by definition we
have $[\Omega] \in F^2_XH^2_{DR}(\wt{X})$.

Applying the restriction map $P:H^2(\wt{X},\C)) \to H^2(Z_x,\C)$, we
deduce that $P([\Omega]) \in F^2H^2(Z_x,\C) \cap
\overline{F}^1(Z_x,\C)$. But since the mixed Hodge structure on
$H^2(Z_x,\C)$ has weights $\leq 2$, this intersection is trivial.
\endproof

\begin{corr}\label{ex}
Let $\pi:\wt{X} \to X$ be a smooth projective resolution of a
symplectic variety $X$. Then for every closed point $x \in X$, the
restriction of the $2$-form $\Omega$ on $\wt{X}$ to the formal
neighborhood of the preimage $\pi^{-1}(x) \subset \wt{X}$ is exact.
\end{corr}

\proof{} We again take $k=\C$ and pass to the analytic topology. It
suffices to find an open neighborhood $U \subset X$ of the point $x$
such that $\Omega$ is exact on $\pi^{-1}(U) \subset \wt{X}$. By
Lemma~\ref{from.b}, the cohomology class $[\Omega]$ of the form
$\Omega$ comes from a cohomology class $[\Omega_X] \in
H^2(X,\C)$. Taking a small enough neighborhood $U \subset X$ of the
point $x \in X$, we can insure that $[\Omega_X] = 0$ in
$H^2(U,\C)$. Thus we can take a neighborhood $U \subset X$ such that
$[\Omega] = 0$ on $\wt{U} = \pi^{-1}(U) \subset X$. Analyzing the
Hodge-de Rham spectral sequence for $\wt{U}$, we see that
$$
\Omega = d_2\beta \mod
d\left(H^0\left(\wt{U},\Omega^1_{\wt{U}}\right)\right),
$$ 
where $d$ is the de Rham differential, $d_2$ is the second
differential in the spectral sequence and $\beta$ is a class in
$H^1(\pi^{-1}(U),\calo_{\wt{U}})$. Since $X$ has rational
singularities, $\beta = 0$ and $\Omega = d\alpha$ for some $1$-form
$\alpha$ on $\wt{U}$.
\endproof

\begin{lemma}\label{pb}
Let $\pi:\wt{X} \to X$ be a smooth projective resolution of
singularities of a symplectic variety $X$. Denote by $\Omega \in
\Omega^2(\wt{X})$ the symplectic form on the manifold $\wt{X}$.  Let
$\sigma:Z \to U$ be a smooth map of smooth connected algebraic
manifolds, and assume given a commutative square
\begin{equation}\label{cd}
\begin{CD}
Z @>{\eta}>> \wt{X}\\
@V{\sigma}VV @VV{\pi}V\\
U @>{\eta_0}>> X.
\end{CD}
\end{equation}
Then there exists a dense open subset $U_0 \subset U$ and $2$-form
$\Omega_U \in \Omega^2(U_0)$ on $U_0$ such that
\begin{equation}\label{rv}
\sigma^*\Omega_U = \eta^*\Omega
\end{equation}
on $\sigma^{-1}(U_0) \subset Z$.
\end{lemma}

\proof{} Consider the fibered product $\overline{Z} = U \times_X
\wt{X}$ equipped with the reduced scheme structure. Since $\wt{X}
\to X$ is projective, $\wt{Z}$ is projective over $U$. Resolving the
singularities, we can find a smooth, possibly non-connected scheme
$Z_0$ projective over $U$ and dominating every irreducible component
of the scheme $\overline{Z}$. Consider the fibered product $Z_0
\times_{\overline{Z}} Z_0$. Resolving its singularities, we can find
a smooth scheme $Z_1$ projective over $U$ which dominates every
irreducible component of $Z_0 \times_{\overline{Z}} Z_0$.  Composing
the map $\Z_1 \to Z_0 \times_{\overline{Z}} Z_0$ with the natural
projections, we obtain two maps $p_1,p_2:Z_1 \to Z_0$. By the
Bertini Theorem, we can replace $U$ with a dense open subset so that
the projection $Z_0 \to U$ and both maps $p_1,p_2:Z_1 \to Z_0$
become smooth.

The schemes $Z_1$, $Z_0$ with the map $p_1,p_2:Z_1 \to Z_0$ are the
first terms of a simplicial hypercovering of the scheme $\wt{Z}/U$,
see \cite{D}. In particular, let $Z_{1,u}$, $Z_{2,u}$
$\overline{Z}_u$ be the fibers of the schemes $Z_1$, $Z_0$, $\wt{Z}$
over a closed point $u \in U$; then for every integer $l \geq 0$ we
have an exact sequence
\begin{equation}\label{spc}
\begin{CD}
0 @>>> \gr_l^WH^l(\wt{Z}_u,\C) @>>> H^l(Z_{0,u},\C) @>{p_1^* -
p_2^*}>> H^l(Z_{1,u},\C),
\end{CD}
\end{equation}
where $\gr_p^WH^p(\wt{Z}_u,\C)$ is the weight-$p$ part of the
cohomology group $H^p(\wt{Z},\C)$ with respect to the weight
filtration associated to the mixed Hodge structure.

Now, the form $\Omega_U$ is uniquely defined by the conditions of
the Lemma. Therefore it suffices to construct it locally on $U$, and
even locally in \'etale topology. Moreover, since $Z$ is smooth, it
suffices to check \eqref{rv} generically on $U$ and on
$Z$. Therefore, possibly replacing $U$ with its \'etale cover and
shrinking $Z$, we can assume that the map $\eta:Z \to \wt{X}$
factors through $Z_0 \to \wt{X}$. Thus it suffices to prove the
Lemma when $Z \subset Z_0$ is a connected component of the scheme
$Z_0$.

To do this, we first note that by Lemma~\ref{from.b}, the de Rham
cohomology class $\eta^*[\Omega]$ of the form $\Omega$ vanishes on
the fibers of the map $Z_0 \to U$. Since $Z_0/U$ is projective, the
Hodge to de Rham spectral sequence for the relative de Rham
cohomology $H^\hdot_{DR}(Z_0/U)$ degenerates; therefore the form
$\eta^*\Omega$ itself vanishes on fibers of the map $Z_0 \to U$. We
conclude that for any closed point $u \in U$ and any tangent vector
$\xi \in T_uU$, the relative $1$-form $\alpha_\xi=\eta^*\Omega
\cntrct \xi \in H^0(Z_{0,u},\Omega^1(Z_{0,u}))$ is well-defined
(locally on $Z_0$, lift $\xi$ to a vector field on $Z_0$, and check
that $\alpha_\xi$ does not depend on the lifting). The same argument
applies to the map $Z_1 \to U$. Moreover, since $\eta \circ p_1 =
\eta \circ p_2:Z_1 \to \wt{X}$, we have
$$
p_1^*\alpha_\xi = p_1^*\eta^*\Omega \cntrct \xi = p_2^*\eta^*\Omega
\cntrct \xi = p_2^*\alpha_\xi.
$$
Applying \eqref{spc} with $l=1$, we conclude that the cohomology
class $[\alpha_\xi] \in H^1_{DR}(Z_{0,u})$ of the form $\alpha_\xi$
comes from a class in $H^1(\overline{Z}_u,\C)$. But we know that
$R^1\pi_*\C = 0$ on $X$, and by proper base change, this implies
$H^1(\overline{Z}_u,\C) = 0$, so that $[\alpha_\xi]=0$. Since the
Hodge to de Rham spectral sequence for $Z_{0,u}$ degenerates, we
conclude that $\alpha_\xi = 0$, for every closed point $u \in U$ and
any tangent vector $\xi \in T_uU$. In other words, for every closed
point $z \in Z \subset Z_0$ we have $\eta^*\Omega(\xi_1,\xi_2)=0$ if
at least one of the tangent vectors $\xi_1,\xi_2 \in T_zZ$ is
vertical with respect to the map $\sigma:Z \to U$. This means that
$$
\eta^*\Omega \in H^0(Z,\sigma^*\Omega^2(U)) \subset
H^0(Z,\Omega^2(Z)).
$$
To finish the proof, note that since $Z$ is connected and projective
over $U$, the natural map $\sigma^*:H^0(U,\Omega^2(U)) \to
H^0(Z,\sigma^*\Omega^2(U))$ is one-to-one.
\endproof

\subsection{Proofs.}
To prove Theorem~\ref{main}, we have to apply the results of the
last Subsection to a particular resolution of singularities of the
normal symplectic variety $X$. We will use the canonical resolution
of singularities $\pi:\wt{X} \to X$ constructed, for example, in
\cite{BM} or in \cite{EH}. It enjoys the following two crucial
properties:
\begin{enumerate}
\item The map $\pi:\wt{X} \to X$ is one-to-one over the smooth part
$X^o \subset X$.
\item Every vector field $\xi$ on $X$ lifts to a vector field
$\wt{\xi}$ on $\wt{X}$.
\end{enumerate}
Actually, the lifting property with respect to the vector fields is
not claimed either in \cite{BM} or in \cite{EH} -- the authors only
prove equivariance with respect to automorphisms. However, the
lifting property is easily deduced from this, for instance, along
the lines of \cite[Lemma 2.2]{K1}.

\proof[Proof of Theorem~\ref{main}.]
Let $X$ be a normal symplectic variety, and let $\pi:\wt{X} \to X$
be the canonical resolution of singularities of the variety
$X$. Since all the claims of the Theorem are local, we may assume
that $X$ is affine.

To prove that $X$ is locally exact, it suffices to apply
Corollary~\ref{ex} to the resolution $\wt{X}$ and then invoke
Lemma~\ref{exl}.

To prove that $X$ is holonomic, we use an idea that essentially goes
back to \cite{hu}. For any locally closed subscheme $U \subset X$,
denote by $\pi^{-1}(U) \subset \wt{X}$ its preimage under the map
$\pi:\wt{X} \to X$ equipped with the reduced scheme structure.  Let
$Y \subset X$ be an integral Poisson subscheme, and let $U \subset
Y$ be the open dense subset such that $U$ is smooth and the map
$\pi:\pi^{-1}(U) \to U$ is generically smooth on $\pi^{-1}(u)
\subset \pi^{-1}(U)$ for any closed point $u \in U$. It suffices to
prove that the Poisson bivector $\Theta$ is non-degenerate on the
cotangent space $T^*_y(Y)$ for every point $y \subset U$. Fix such a
point $y$. If $\Theta$ is degenerate, then some non-trivial covector
in $T^*_y(Y)$ lies in its annihilator. In other words, there exists
a function $f$ on $X$ such that $df \neq 0$ at $y \in U$, but the
Hamiltonian vector field
$$
H_f = df \cntrct \Theta
$$
vanishes at the point $y \in Y \subset X$. On the other hand,
generically on $X$ we have a well-defined symplectic form $\Omega$,
and we have
$$
df = H_f \cntrct \Omega.
$$
Since the resolution $\wt{X}$ is canonical, the vector field $H_f$
lifts to a vector field $\wt{H}_f$ on $\wt{X}$. Generically we have
$$
\pi^*(df) = \wt{H}_f \cntrct \Omega.
$$
But $\wt{X}$ is smooth and $\Omega$ is defined everywhere. Therefore
this equality also holds everywhere -- in particular, on every
connected component of the smooth part of $\pi^{-1}(U)$. By
Lemma~\ref{pb}, we can replace $U$ with a dense open subset so that
on such a connected component we have $\Omega = \pi^*\Omega_0$ for
some $2$-form $\Omega_0$ on $U$. Since $H_f$ preserves $U \subset
X$, the vector field $\wt{H}_f$ preserves (that is, is tangent to)
its preimage $\pi^{-1}(U)$, and we have
$$
\wt{H}_f \cntrct \Omega = H_f \cntrct \Omega_0.
$$
Since $H_f$ vanishes at $y$, we conclude that $df = 0$ in every
point in the smooth part of $\pi^{-1}(y) \subset X$. By assumption
$df \neq 0$ at $y \in U$, and this contradicts our smoothness
assumptions on the map $\pi:\pi^{-1}(U) \to U$.

Finally, we have to prove that the normalization $Y^{nrm}$ of an
irreducible Poisson subscheme $Y$ is a symplectic variety. We note
that since we already know that $X$ is holonomic, $Y$ is generically
symplectic. Moreover, every Poisson subscheme in $Y$ is also a
Poisson subscheme in $X$. Therefore $Y$ itself is holonomic, and so
is its normalization $Y^{nrm}$. By Lemma~\ref{sm} this means that
the regular part $Y^{reg} \subset Y^{nrm}$ carries a symplectic
form. We have to prove that this form extends to a smooth resolution
of the variety $Y^{nrm}$.

To do this, note that generically on $Y$ and locally in \'etale
topology, the projection $\pi^{-1}(Y) \to Y$ admits a section. More
precisely, taking a sufficiently small smooth open dense subset $U
\subset Y$, we can assume that there exists a Galois cover
$\kappa:U' \to U$ and a map $\sigma:U' \to \wt{X}$ such that $\sigma
\circ \pi = \kappa:U' \to U \subset X$. By Lemma~\ref{pb}, shrinking
$U$ even further we can assume that $\sigma^*\Omega =
\kappa^*\Omega_U$ for some $2$-form $\Omega_U$ on $U$. We will prove
that (1) this form $\Omega_U$ extends to a resolution $\wt{Y}$ of
the variety $Y$, and (2) at least generically on $Y$, the form
$\Omega_U$ coincides with the form given by the Poisson structure on
$Y$.

\smallskip

\proof[Step {\normalfont 1}: $\Omega_U$ extends to a form on a
resolution $\wt{Y} \to Y$.]

Let $Y'$ be the normalization of the scheme $Y$ in the Galois cover
$U' \to U$, and consider the fibered product $\wt{X} \times_Y Y'$
equipped with the reduced scheme structure. Let
$\overline{\sigma(U')} \subset \wt{X} \times_Y Y'$ be the closure of
the image of the section $\sigma:U' \to \wt{X} \times_Y Y'$, and let
$\wt{Y'}$ be a smooth projective $Gal(U'/U)$-equivariant resolution
of singularities of the closure $\overline{\sigma(U')}$.

We have a $Gal(U'/U)$-equivariant smooth resolution $\wt{Y'} \to Y'$
and a map $\sigma:\wt{Y'} \to \wt{X}$.  By assumption, over $U'
\subset Y'$ the $2$-form $\sigma^*\Omega$ on $\wt{Y'}$ coincides
with $\kappa^*\Omega_U$. In particular, the $2$-form
$\sigma^*\Omega$ is $Gal(U'/U)$-invariant. The quotient
$\wt{Y}_0=\wt{Y'}/Gal(U'/U)$ is a normal algebraic variety equipped
with a projective birational map onto $Y$.

In general, let $f:Z' \to Z$ be a finite morphism between normal
algebraic varieties such that $Z'$ is smooth and equipped with an
action of a finite group $G$, and $f:Z' \to Z$ is generically a
Galois cover with Galois group $G$. Then for any $p \geq 0$, any
$G$-invariant $p$-form $\alpha$ on $Z'$ gives by descent a $p$-form
on an open smooth subset of the variety $Z$, and this form extends
to any smooth projective resolution $\wt{Z} \to Z$ (this is
well-known; for a sketch of a proof see e.g. \cite[Lemma
3.3]{K3}). In particular, take $p=2$, $Z' = \wt{Y'}$ and $G =
Gal(U'/U)$. Then the $Gal(U'/U)$-invariant $2$-form $\sigma^*\Omega$
on $\wt{Y'}$ gives a $2$-form $\Omega$ on any smooth projective
resolution $\wt{Y} \to \wt{Y}_0$.

Since generically on $\wt{Y'}$ we have $\sigma^*\Omega =
\kappa^*\Omega_U$, we conclude that the $2$-form $\Omega_U$ on $U
\subset Y$ extends to a smooth projective resolution $\wt{Y} \to
\wt{Y}_0 \to Y$.

\smallskip

\proof[Step {\normalfont 2}: $\Omega_U$ is compatible with the
Poisson structure on $Y$.]

Since the Poisson structure on $U$ is non-degenerate, the tangent
bundle $\T(U)$ is generated by Hamiltonian vector fields. Since the
map $\kappa:U' \to U$ is \'etale, all vector fields on $U$ lift to
vector fields on $U'$. Thus it suffices to check that for every two
Hamiltonian vector fields
$$
H_f = df \cntrct \Theta, \qquad H_g = dg \cntrct \Theta,
$$
we have $(\sigma^*\Omega)(H_f,H_g) = \{f,g\}$ on $U'$. Again, both
$H_f$ and $H_g$ lift to vector fields $\wt{H}_f$, $\wt{H}_g$ on
$\wt{X}$, and by Lemma~\ref{pb} we have
$$
\left(\sigma^*\Omega\right)(H_f,H_g) =
\sigma^*\left(\Omega\left(\wt{H}_f,\wt{H}_g\right)\right).
$$
Thus it suffices to check that
$$
\Omega\left(\wt{H}_f,\wt{H}_g\right) = \pi^*\{f,g\}.
$$
But this equation makes sense everywhere on $\wt{X}$. Since $\wt{X}$
is reduced, it suffices to check it generically, where it follows
from the definition of the form $\Omega$.
\endproof

\subsection{Symplectic resolutions.}\label{crep}

Assume now given a symplectic variety $X$ and a smooth projective
resolution $\wt{X} \to X$ which is {\em crepant} -- in this context,
it means that the canonical bundle $K_{\wt{X}}$ is trivial. Since
the top degree $\Omega^{\frac{1}{2}\dim X}$ of the symplectic form
is a section of $K_{\wt{X}}$, the form $\Omega_{\wt{X}}$ in this
case must be non-degenerate everywhere, not only at the generic
point of the scheme $\wt{X}$. It turns out that this imposes strong
restrictions on the geometry of $\wt{X}$. We start with the
following general fact.

\begin{lemma}\label{vnsh}
Let $X$, $Y$ be algebraic varieties over $k$, assume that $X$ is
smooth, and let $\tau:X \to Y$ be a projective map. Then
$R^p\tau_*\Omega^q_X = 0$ whenever $p+q > X \times_Y X$.
\end{lemma}

\proof{} (I am grateful to H. Esnault and E. Viehweg for suggesting
and explaining this proof to me.)

The claim is local in $Y$, so that we may assume that $Y$ is
affine. Choose a projective variety $\overline{Y}$, a smooth
projective variety $\overline{X}$ and a projective map
$\tau:\overline{X} \to \overline{Y}$ so that $Y \subset
\overline{Y}$ is a dense open subset in $\overline{Y}$, $X =
\pi^{-1}(Y) \subset \overline{X}$ is a dense open subset in
$\overline{X}$, $\overline{\tau}:\overline{X} \to \overline{Y}$
extends the given map $\tau:X \to Y$, and the complement $E =
\overline{X} \setminus X$ is a simple normal crossing divisor in
$\overline{X}$. Choose an ample line bundle $M$ on
$\overline{Y}$. Consider the logarithmic de Rham complex
$\Omega^\hdot_{\overline{X}}\langle \log E \rangle$. Let $l >\!\!>
0$ be an integer large enough so that the sheaves
\begin{equation}\label{dr.im}
R^p\overline{\tau}_*\Omega^q_{\overline{X}}\langle \log E
\rangle(-E) \otimes M^{\otimes l}
\end{equation}
are acyclic and globally generated for all $p,q$. Replace $M$ with
$M^{\otimes l}$. Then by \cite[Corollary 6.7]{EV}, we have
$$
H^p\left(\overline{X},\Omega^q_{\overline{X}}\langle \log E \rangle
\otimes \tau^*M^{-1}\right) = 0
$$
whenever $p+q < \dim X - r(\tau)$. Here $r(\tau)$ is a certain
constant defined in \cite[Definition 4.10]{EV}; although the
definition is textually different, it is immediately obvious that
$\dim X + r(\tau) = \dim X \times_Y X$. By Serre duality (see
\cite[Remark 6.8 b)]{EV}),
$$
H^p\left(\overline{X},\Omega^q_{\overline{X}}\langle \log E
\rangle(-E) \otimes \tau^*M\right) = 0
$$
whenever $p+q > \dim X + r(\tau) = \dim X \times_Y X$. Since the
sheaves \eqref{dr.im} are acyclic, the Leray spectral sequence for
these cohomology groups collapses, and we conclude that
$$
H^0\left(R^p\overline{\tau}_*\Omega^q_{\overline{X}}\langle \log
E \rangle(-E) \otimes M\right) = 0
$$
wherever $p+q > \dim X \times_Y X$. Since the sheaves \eqref{dr.im}
are globally generated, this implies
$R^p\overline{\tau}_*\Omega^q_{\overline{X}} \langle \log E
\rangle(-E)=0$. Restricting to $Y \subset \overline{Y}$, we get the
claim.
\endproof

This can applied to the symplectic situation because of the
following.

\begin{lemma}\label{semismall}
Let $\pi:\wt{X} \to X$ be a projective birational map from a smooth
variety $\wt{X}$ with $K_{\wt{X}}=0$ to a symplectic variety
$X$. Then the map $\pi:\wt{X} \to X$ is semismall, in other words,
$\dim \wt{X} \times_X \wt{X} = \dim X$.
\end{lemma}

\proof{} For any $p \geq 0$, let $X_p \subset X$ be the closed
subvariety of points $x \in X$ such that $\dim \pi^{-1}(x) \geq
p$. It suffices to prove that $\codim X_p \geq 2p$. By
Lemma~\ref{pb}, there exists an open dense subset $U \subset X_p$
such that the restriction $\Omega_F$ of the form $\Omega =
\Omega_{\wt{X}}$ onto every connected component $F$ of the smooth
part of the set-theoretic preimage $\pi^{-1}(U) \subset \wt{X}$
satisfies
$$
\Omega_F = \pi^*\Omega_U
$$
for some $2$-form $\Omega_U \in H^0(U,\Omega^2_U)$. Therefore the
rank $\rk \Omega_F$ satisfies $\rk \Omega_F \leq \dim U$. On the
other hand, since $K_{\overline{X}}$ is a trivial line bundle, the
form $\Omega$ is non-degenerate on $\wt{X}$, and we have $\rk
\Omega_F \geq \dim F - \codim F$. By definition of $X_p \subset X$,
we can choose a component $F$ such that $\dim F = \dim U + p$. Then
together these two inequalities give $\codim X_p = \codim F + p \geq
\dim F - \dim U + p = 2p$, as required.
\endproof

\begin{theorem}\label{sympl}
Let $\pi:\wt{X} \to X$ be a projective birational map with smooth and
symplectic $\wt{X}$. Let $x \in X$ be a closed point, and let $E_x =
\pi^{-1}(x) \subset \wt{X}$ be the set-theoretic fiber over the point
$x$. Then for odd $k$ we have $H^k(E_x,\C)=0$, while for even $k=2p$
the Hodge structure on $H^k(E_x,\C)$ is pure of weight $k$ and Hodge
type $(p,p)$.
\end{theorem}

\begin{lemma}\label{hdg}
Let $p$ be an integer, and let $V$ be an $\R$-mixed Hodge structure
with Hodge filtration $F^\hdot$ and weight filtration
$W_\idot$. Assume that $W_{2p}V = V$ and $F^pV=V$. Then $V$ is a
pure Hodge-Tate structure of weight $2p$ (in other words, every
vector $v \in V$ is of Hodge type $(p,p)$).
\end{lemma}

\proof{} Since $V = F^pV$, the same is true for all associated
graded pieces of the weight filtration on $V$. Therefore we may
assume that $V$ is pure of weight $k \leq 2p$. If $k < 2p$, we must
have $V = F^pV \cap \overline{F^pV} = 0$, which implies $V = 0$. If
$k = 2p$, the same equality gives $V = V^{p,p}$.
\endproof

\proof[Proof of Theorem~\ref{sympl}.] By Lemma~\ref{semismall},
Lemma~\ref{vnsh} applies to $\pi:\wt{X} \to X$ and shows that
\begin{equation}\label{est}
R^p\pi_*\Omega^q_{\wt{X}} = 0
\end{equation}
whenever $p + q > \dim \wt{X}$. Since $\wt{X}$ is symplectic, we have an
isomorphism $\T_{\wt{X}} \cong \Omega^1_{\wt{X}}$ between the
tangent and the cotangent bundle on $\wt{X}$. This implies that
$\Omega^q_{\wt{X}} \cong \Omega^{\dim \wt{X} - q}_{\wt{X}}$, and
\eqref{est} also holds whenever $p > q$.

Denote by $\X$ the completion of the variety $\wt{X}$ in the closed
subscheme $E_x = \pi^{-1}(x) \subset \wt{X}$. Since the map
$\pi:\wt{X} \to X$ is proper, by proper base change the group
$$
H^p(\X,\Omega^q_\X)
$$
for any $p$, $q$ coincides with the completion of the stalk of the
sheaf $R^p\pi_*\Omega_{\overline{X}}^q$ at the point $x \in X$. Therefore
$H^p(\X,\Omega^q_\X) = 0$ whenever $p > q$. The stupid filtration on
the de Rham complex $\Omega^\hdot_\X$ of the formal scheme $\X$
induces a descreasing filtration $F^\hdot$ on the de Rham cohomology
groups $H_{DR}(\X)$ which we call the weak Hodge filtration. Of
course, the associated spectral sequence does not
degenerate. Nevertheless, since $H^p(\X,\Omega^q_\X) = 0$ when $p >
q$, we have $H^k_{DR}(\X) = F^pH^k_{DR}(\X)$ whenever $k \leq 2p$.

It is well-known that the canonical restriction map
$$
H^\hdot_{DR}(\X) \to H^\hdot(E_x,\C)
$$
is an isomorphism. By definition (see \cite{D}), to obtain the Hodge
filtration on the cohomology groups $H^\hdot(E_x)$, one has to
choose a smooth simplicial resolution $E^\hdot_x$ for the variety
$E_x$ and take the usual Hodge filtration on
$H^\hdot{DR}(E^\hdot_x)$. The embedding $E_x \to \X$ gives a map
$E^\hdot_x \to \X^\hdot$, where $\X^\hdot$ is $\X$ considered as a
constant simplicial variety. The corresponding restriction map
$$
H^\hdot_{DR}(\X) \to H^\hdot_{DR}(E^\hdot_x)
$$
is also an isomorphism, and it sends the weak Hodge filtration on
the left-hand side into the usual Hodge filtration on the right-hand
side. We conclude that $H^k(E_x) = F^p(E_x)$ whenever $k \leq
2p$. It remains to recall that by definition, we have $H^k(E_x) =
W_kH^k(E_x)$, and apply Lemma~\ref{hdg}.
\endproof

To conclude this subsection, we would like to note that in the
particular case when $X = T^*(G/B)$ is the Springer resolution of
the nilpotent cone $Y = \N \subset \g^*$ in the coadjoint
representation $\g^*$ of a semisimple algebraic group $G$,
Theorem~\ref{sympl} has been already proved by C. de Concini,
G. Lusztig and C. Procesi in \cite{lu}. They proceed by a direct
geometric argument. As a result, they obtain more: not only do the
cohomology groups carry a Hodge structure of Hodge-Tate type, but in
fact they are spanned by cohomology classes of algebraic
cycles. This is true even for cohomology groups with integer
coefficients. Motivated by this, we propose the following.

\begin{conj}
In the assumptions of Theorem~\ref{sympl}, the cohomology groups
$H^k(E_x,\Z)$ are trivial for odd $k$, and are spanned by cohomology
classes of algebraic cycles for even $k$.
\end{conj}

We also expect that an analogous statement holds for $l$-adic
cohomology groups, possibly even over fields of positive
characteristic.

\section{Stratification and product decomposition.}\label{str}

We now turn to the algebraic study of Poisson schemes.

\begin{prop}\label{holo}
Let $X$ be a Noetherian integral Poisson scheme over $k$. Assume
that $X$ is excellent as a scheme and holonomic as a Poisson
scheme. Then there exists a stratification $X_i \subset X$ by
Poisson subschemes such that the open parts $X_i^o$ of the strata
are smooth and symplectic. The stratification is canonical -- in
particular, it is preserved by all automorphisms of the scheme $X$
and by all vector fields. The only integral Poisson subschemes in
$X$ are the irreducible components of the closed strata $X_i$.
\end{prop}

\proof{} Let $Y \subset X$ be singular locus of the scheme
$X$. Since $X$ is excellent, $Y$ is a proper closed subscheme
preserved by all automorphisms of $X$ and by all vector fields. In
particular, it is presevred by all Hamiltonian vector fields. This
means that $Y \subset X$ is a Poisson subscheme. It is automatically
holonomic. Since $\dim Y < \dim X$, we can stratify it by
induction. Thus it suffices to stratify the smooth part $X \setminus
Y \subset X$.  In this case all the claims follow from
Lemma~\ref{sm}.
\endproof

Note that this immediately implies that if a holonomic Poisson
scheme $X$ is locally exact, then every Poisson subscheme $Y \subset
X$ is also locally exact. Indeed, every such subscheme must be a
closed stratum $X_i$, hence it is preserved by {\em all} locally
defined vector fields on $X$ -- in particular, by a vector field
$\xi$ satisfying \eqref{dilat}. 

Next, we construct the product decomposition \eqref{prod}. We need
the following general result.

\begin{lemma}\label{triv}
Let $M$ be a vector space equipped with a descreasing filtration
$F^pM$, $p \geq 0$ such that $\codim F^pM < \infty$ and $M$ is
complete with respect to the topology induced by $F^\hdot$. Assume
that $M$ is a module over the algebra $A = k[[x_1,\dots,x_n]]$ of
formal power series in $n$ variables, and that $x_i \cdot F^pM
\subset F^{p+1}M$ for every $i$, $p$. Finally, assume that the
module $M$ is equipped with a flat connection $\nabla:M \to M
\otimes \Omega^1_A$, and let $M^\nabla \subset M$ be the subspace of
flat sections. Then the natural map
$$
M^\nabla \whotimes A \to M
$$
from the completed tensor product $M \whotimes A$ to $M$ is an
isomorphism.
\end{lemma}

\proof{} This is completely standard, but under an additional
assumption that $M$ is finitely generated. We give a proof to show
that our assumptions are in fact sufficient.

Consider the algebra $D$ of differential operators $A \to A$; as a
vector space, $D \cong A[\xi_1,\dots,\xi_n]$, where $\xi_i$ denotes
the differential operator $\frac{\6}{\6x_i}$. An $A$-module equipped
with a flat connection $\nabla$ is the same as a left $D$-module
(the generator $\xi$ acts by covariant derivative with respect to
$\frac{\6}{\6x_i}$). Since $M$ is a cocompact topological vector
space, we have $M=N^*$, where $N$ is the (discrete) vector space of
continous linear maps $M \to k$. The filtration $F$ induces an
increasing filtation $F_pN$, $p \geq 0$ such that $\dim F_p N <
\infty$. The maps $x_\idot,\xi_\idot:M \to M$ induce by duality maps
$x_\idot,\xi_\idot:N \to N$ satisfying the same commutation
relations. We have $x_i \cdot F_{p+1}N \subset F_pN$ for every $i$,
$p$; therefore every element $a \in N$ is annihilated by a high
power of every $x_i$, and $N$ becomes a left $D$-module.

Let $N^o \subset N$ be the common kernel of multiplication by
$x_1,\dots,x_n$. We first prove that for any filtered left
$D$-module $N$ satisfying $x_i \cdot F_{p+1}N \subset F_pN$, the
natural map
$$
a_{N,N^o}:N^o \otimes k[\xi_1,\dots,\xi_n] \to N
$$
induced by the $k[\xi_1,\dots,\xi_n]$-module structure on $N$ is an
isomorphism.

By induction, it suffices to consider the case $n=1$. Indeed, let
$N' \subset N$ be the kernel of multiplication by $x_1$. Then $N'$
carries a natural structure of a filtered module over
$k[[x_2,\dots,x_n]]$ equipped with a flat connection, and it
satisfies all our assumptions. The map $a_{N,N^o}$ factors as
$$
\begin{CD}
N^o \otimes k[\xi_1,\dots,\xi_n] @>{a_{N',N^o} \otimes \id}>> N'
\otimes k[\xi_1] @>{a_{N,N'}}>> N,
\end{CD}
$$
we know by the inductive assumption that $a_{N',N^o}$ is an
isomorphism, and we have to prove that $a_{N,N'}$ is also an
isomorphism. Thus we may forget about $x_i$, $\xi_i$ for $i \geq 2$
and assume $n=1$, $N^o = \Ker x_1 \subset N$.

To simplify notation, let $a_N = a_{N,N^o}$.  Note that for every
$k[[x_1]]$-module $N$ satisfying our assumptions, the kernel $\Ker
x_1$ must be non-trivial (for instance, it contains the smallest
non-trivial term in the filtration $F_\idot N$). In particular, this
applies to the kernel $\Ker a_N \subset N^o \otimes k[\xi_1]$ of the
map $a_N$. But the kernel of $x_1$ on $\Ker a_N \subset N^o \otimes
k[\xi_1]$ coincides with $\Ker a_N \cap \Ker x_1 \subset N^o \otimes
k[\xi_1]$, which in turn is equal to the kernel of the map $a_N$ on
$N^o = \Ker x_1 \subset N \otimes k[\xi_1]$; since by definition
this kernel is trivial, we must also have $\Ker a_N = 0$. We
conclude that $a_N$ is injective. It remains to check that $N$ is
generated by $\Ker x_1$ as a $k[\xi_1]$-module. Denote $N_p = \Ker
x_1^p \subset N$. By assumption $F_pN \subset N_p$, so that $N =
\bigcup N_p$. Thus by induction is suffices to prove that
$$
N_p \subset N_{p-1} + k[\xi_1] \cdot N_1.
$$
This is immediate. For every element $m \in N_p$, let
$$
m_0 = m - \frac{1}{(p-1)!}\xi_1^{p-1}x_1^{p-1}m.
$$
Since $x_1\xi_1 - \xi_1x_1 = \id$, we have $x_1^{p-1}m_0 = 0$, and
$m_0 \in N_{p-1}$.

Now, we have proved that $N \cong V \otimes k[\xi_1,\dots,\xi_n]$
for some vector space $V=N^o$; since $x_i$ vanishes on $N^o$, both
the operators $x_\idot$ and the operators $\xi_\idot$ act on $V
\otimes k[\xi_1,\dots,\xi_n]$ via the second factor. Therefore
$$
M = N^* \cong V^* \whotimes A,
$$
where again $x_\idot$ and $\xi_\idot$ act on the product via the
second factor. To prove the Lemma, it suffices to show that the
natural map $M^\nabla \to M$ identifies $M^\nabla$ with $V^* \otimes
1 \subset V^* \whotimes A$. If $V^*$ is one-dimensional, this is
obvious: $A^\nabla$ is indeed the line spanned by $1 \in A$. But the
completed tensor product functor $W \mapsto W \whotimes A$ is exact
and commutes with arbitrary inverse limits, and in particular, with
arbitrary products; the flat sections functor $M \mapsto
M^\nabla$ is left-exact and therefore also commutes with arbitrary
products. Since every cocompact vector space $V^*$ is a (possibly
infinite) product of one-dimensional vector spaces, we are done.
\endproof

\begin{prop}\label{prodo}
Let $A$ be a complete Poisson local algebra over $k$ with maximal
ideal $\m \subset A$, and assume given a prime Poisson ideal $J
\subset A$ such that the quotient $A/J$ is a regular complete local
algebra with with non-degenerate Poisson structure. Then there
exists a complete local Poisson algebra $B$ and a Poisson
isomosphism 
\begin{equation}\label{whpro}
A \cong B \whotimes_k (A/J)
\end{equation}
between the algebra $A$ and the completed tensor product of the
algebras $A/J$ and $B$. Moreover, the Poisson scheme $\Spec A$ is
holonomic if and only if the Poisson scheme $\Spec B$ is holonomic,
and the Poisson algebra $A$ is exact if and only if the Poisson
algebra $B$ is exact.
\end{prop}

\proof{} For every integer $d \geq 1$, denote by $W_d$ the power
series algebra $k[[x_1,\dots,x_k,y_1,\dots,y_d]]$ with the standard
Poisson structure induced by the symplecic form $dx_1 \wedge dy_1 +
\dots + dx_d \wedge dy_d$. We first prove that there exists a
Poisson isomorphism
\begin{equation}\label{dim1}
A \cong W_1 \whotimes A'
\end{equation}
for some complete local Poisson algebra $A'$. By the formal Darboux
Theorem, there exists a Poisson isomorphism $A/J \cong W_d$, where
$2d = \dim A/J$. Fix arbitrary liftings $f,g \in A$ of $x_1,y_1 \in
A/J$. Then $\{f,g\} = 1 \mod J$. We claim that there exist a series
of functions $f_l \in J^l$, $l \geq 1$ such that
\begin{equation}\label{g_l}
\{f,g+g_1+\dots+g_l\} = 0 \mod J^{l+1}.
\end{equation}
Indeed, let $\xi(a) = \{a,g\}$ for $a \in A$; then by induction on
$l$ it suffices to prove that $\xi:J^l/J^{l+1} \to J^l/J^{l+1}$ is
surjective. However, on $J^l/J^{l+1}$ we have $\xi(fa)=f\xi(a)+a$,
so that $\xi$ induces a flat connection on $J^l/J^{l+1}$ considered
as a $k[[f]]$-module. Therefore we can equip $J^l/J^{l+1}$ with
$\m$-adic filtration and apply Lemma~\ref{triv}.

Having chosen a sequence $g_l$ as in \eqref{g_l}, replace $g$ with
$$
g' = g + g_1 + \dots,
$$
so that $\{f,g'\} = 1$ in the algebra $A$, and let
$\xi_1(a)=\{f,a\}$, $\xi_2(a)=\{g',a\}$ for any $a \in A$. Then we
have $k[[f,g']] \cong W_1$, and $\xi_1$ and $\xi_2$ induce a flat
connection on $A$ considered as a $k[[f,g']]$-module. Let $A' = A_0
= \Ker \xi_1 \cap \Ker \xi_2$ and apply Lemma~\ref{triv} to $A$
equipped with $\m$-adic filtration. The space $A'$ is obviously a
complete Poisson algebra, and by Lemma~\ref{triv} we indeed have the
isomorphism \eqref{dim1}. To finish the proof of the first claim,
apply induction on $2d = \dim A/J$.

To prove the second claim, note that by Lemma~\ref{triv} every
Poisson ideal $I \subset A$ is equal to $I_0 \otimes A/J$ for some
Poisson ideal $I_0 \subset B$. Thus every Poisson subscheme in $Y
\subset \Spec A$ is the (completed) product $Y_0 \times \Spec A/J$
of a Poisson subscheme $Y_0 \subset \Spec B$ with the non-degenerate
smooth Poisson scheme $\Spec A/J$. This implies that $\Spec A$ is
holonomic if and only if $\Spec B$ is holonomic.

Finally, the algebra $A/J \cong W_d$ is obviously exact -- for
instance, the Euler vector field
$$
\xi_e = \frac{1}{2}\frac{\6}{\6 x_1} + \dots +
\frac{1}{2}\frac{\6}{\6 x_d} + \frac{1}{2}\frac{\6}{\6 y_1} + \dots
+ \frac{1}{2}\frac{\6}{\6 y_d}
$$
satisfies \eqref{dilat}. Thus if $B$ is exact, with a derivation
$\xi_0:B \to B$ satisfying \eqref{dilat}, then the derivation 
$$
\xi = \xi_0 \otimes \id + \id \otimes \xi_e
$$
of the algebra $A \cong B \otimes (A/J)$ also satisfies
\eqref{dilat}, and $A$ is exact. Conversely, assume that we are
given a derivation $\xi:A \to A$ satisfying \eqref{dilat}. The
product decomposition \eqref{whpro} gives in particular a canonical
embedding $A/J \subset A$, thus a direct some decomposition $A \cong
J \oplus (A/J)$. The restriction $\overline{\xi}:A/J \to A$ of the
map $\xi$ to $A/J \subset A$ decomposes as
$$
\overline{\xi} = \xi' + \xi_1,
$$
where $\xi_1:A/J \to A/J$ is a derivation satisfying \eqref{dilat},
and $\xi':A/J \to J$ is a derivation which is also a derivation with
respect to the Poisson bracket,
\begin{equation}\label{poider}
\xi'(\{a,b\}) = \{\xi(a),b\} + \{a,\xi(b)\}
\end{equation}
for any $a,b \in A/J \subset A$. We claim that $\xi' = H_f$ for some
$f \in J$. Indeed, the
Poisson embedding $A/J \cong W_d \subset A$ induces a $W_d$-module
structure on $J \subset A$ and a flat connection $\nabla$ on the
$W_d$-module $J$. Since the Poisson bivector $\Theta$ on $A/J \cong
W_d$ is non-degenerate, we have
$$
\xi' = \alpha \cntrct \Theta
$$
for some $1$-form $\alpha \in J \otimes_{W_d} \Omega^1(W_d/k)$. Then
the equality \eqref{poider} means exactly that the form $\alpha$ is
closed with respect to the connection $\nabla$, $\nabla\alpha=0$. By
Lemma~\ref{triv}, all the higher de Rham cohomology groups of the
flat module $J$ are trivial, so that $\nabla\alpha=0$ implies that
$\alpha = \nabla f$ for some $f \in J$. Thus in turn means that we
indeed have
$$
\xi' = df \cntrct \Theta = H_f.
$$
Replacing $\xi:A \to A$ with $\xi - H_f:A \to A$, we obtain a
derivation that still satisfies \eqref{dilat}, but now also
preserves $A/J \subset A$. Therefore it also preserves $B \subset
A$, and induces a derivation on $B$ satisfying \eqref{dilat}. This
means that $B$ is exact.
\endproof

\begin{remark}
This Proposition is well-known in the theory of Poisson structures
on $C^\infty$-manifolds, see \cite{wein}; the decomposition in this
case is called the {\em Weinstein decomposition}. Our proof is
essentially the same as Weinstein's, but it is re-set in the
algebraic language and works for singular varieties, too.
\end{remark}

\proof[Proof of Theorem~\ref{th.1}.] Almost all the claims follow
immediately from Proposition~\ref{holo} and
Proposition~\ref{prodo}. In Proposition~\ref{prodo} we take $A =
\wh{\calo}_{X,x}$, the algebra of formal germs on functions on $X$
near $x \in X$. The transversal slice $\Y_x$ is the spectrum of the
algebra $B$ provided by Proposition~\ref{prodo}. To prove that it is
a symplectic variety, let $\wh{Y}$ be a resolution of $\Y_x$, and
consider the product $\wt{Y} \times \wh{X^o_i}_x$ as a resolution of
$\wh{X}_x$. Then, since the product decompostion \eqref{prod} is
Poisson, the symplectic form $\Omega$ on this product satisfies
$$
\Omega = p_1^*\Omega_Y + p_2^*\Omega_{X_i},
$$
where $p_1:\wt{Y} \times \wh{X^o_i}_x \to \Y_x$, $p_2:\wt{Y} \times
\wh{X^o_i}_x \to \wh{X^o_i}_x$ are the natural projections, and
$\Omega_Y$, $\Omega_{X_i}$ are the symplectic forms on $\Y_x$ and
$\wh{X^o_i}_x$. Since the forms $\Omega$ and $\Omega_{X_i}$ have no
poles, the form $\Omega_Y$ has no poles either.
\endproof

\section{Group actions.}\label{grp}

We now turn to the proof of Theorem~\ref{th.2}. By
Theorem~\ref{th.1} and Theorem~\ref{main}, we may assume that we are
in the following situation:
\begin{itemize}
\item We have $\Y_x = \Spec A$, where $A$ is a complete local
Poisson algebra, whose maximal ideal $\m \subset A$ is preserved by
all derivations of the algebra $A$. Moreover, the Poisson algebra
$A$ is exact.\label{assu}
\end{itemize}
We will say that a derivation $\xi$ of a Poisson algebra $B$ is {\em
dilating with constant $\theta$} if it satisfies
\begin{equation}
\xi(\{a,b\}) = \{\xi(a),b\} + \{a,\xi(b)\} + \theta\{a,b\}
\end{equation}
for every $a,b\in B$. Since our algebra $A$ is exact, there exist
derivations $\xi:A \to A$ dilating with constant $1$. Fix such a
derivation. As explained in the proof of Lemma~\ref{exl}, to prove
Theorem~\ref{th.2} we essentially need to find a derivation which is
(1) dilating with non-zero constant and (2) can be integrated to an
action of $\gm$ on the algebra $A$. This we will restate and prove
in the following equivalent algebraic form.

\begin{prop}\label{gr}
In the assumptions of \thebull{} above, there exists an
integer $l \neq 0$ and a multiplicative grading
$$
A = \bigoplus_p A^p
$$
on the algebra $A$ such that
$$
\left\{A^p,A^q\right\} \subset A^{p+q-l}.
$$
\end{prop}

\begin{remark}
We may assume that $l$ is positive with any loss of generality. It
would be highly desirable to show that one can choose a grading in
Proposition~\ref{gr} which only has positive weights, $A^p = 0$ for
$p < 0$. Unfortunately, we were unable to prove it -- it seems that
this would require a radically different approach (perhaps a study
of generalized contact singularities would help, see
Remark~\ref{contact}). Thus, Proposition~\ref{gr} is of only limited
use in geometric application.
\end{remark}

\proof{} By \thebull{}, all the derivations of the algebra $A$
preserve the maximal ideal $\m \subset A$ and all its powers $\m^q
\subset A$. In particular, all the ideals $\m^p \subset A$ are
Poisson. Thus for every integer $q \geq 1$ the Artin algebra $A_q =
A/\m^q$ is a Poisson algebra. We have two canonical elements
$$
m,p \in A_q \otimes_k A_q \otimes_k A_q^*
$$ 
-- namely, $m$ defines the multiplication in $A_q$, and $p$ defines the
Poisson bracket.

\begin{lemma}
Let $B$ be a finite-dimensional Poisson algebra over $k$ with
$$
m,p \in B \otimes_k B \otimes_k B^*
$$
giving the multiplication and the Poisson bracket. An endomorphism
$\xi \in \End(B)$ of the vector space $B$ is a dilating derivation
of the algebra $B$ with constant $\theta$ if and only if
$$
\xi(m) = 0 \qquad \xi(p) = \theta p,
$$
where $\xi$ acts on $B \otimes B \otimes B^*$ via the canonical
representation of the Lie algebra $\End_k(B)$.
\end{lemma}

\proof{} Clear. \endproof

\begin{lemma}\label{dil}
Let $\xi \in End(B)$ be a dilating derivation of a
fi\-ni\-te-di\-men\-si\-onal Poisson algebra $B$ with constant
$\theta$. Let $\xi = \xi_s + \xi_n$ be its Jordan decomposition into
the semisimple part $\xi_s$ and the nilpotent part $\xi_n$.

Then both $\xi_s$ and $\xi_n$ are also derivations. Moreover,
$\xi_s$ is dilating with constant $\theta$, and $\xi_n$ is dilating
with constant $0$ (in other words, preserves the Poisson
structure). In particular, $\xi_s \neq 0$.
\end{lemma}

\proof{} By the standard Lie algebra theory, the Jordan
decomposition $\xi = \xi_s + \xi_n$ is universal -- namely, it
induces the Jordan decompositon of the endomorphism $\ad(\xi):V \to
V$ for any finite-dimensional representation $V$ of the reductive
Lie algebra $\End_k(B)$. In particular this applies to the
representation $B \otimes B \otimes B^*$. Now, by assumption $m$ and
$p$ are both eigenvectors of the endomorphism $\ad(\xi)$, with
eigenvalues respectively $0$ and $\theta$. Therefore $\ad(\xi_n)(m)
= \ad(\xi_n)(p) = 0$, and both $m$ and $p$ are eigenvalues of
$\ad(\xi_s)$ with eigenvalues respectively $0$ and $\theta$.
\endproof

Now, by \thebull{} the fixed dilating derivation $\xi:A \to A$
preserves the ideal $\m^q \subset A$ for any $q \geq 1$, so that we
have a dilating derivation $\xi$ of every quotient $A_q = A/\m^q$.

Denote by $T_q \subset GL(A_q)$ the minimal algebraic subgroup whose
Lie algebra $Lie(T_q) \subset \End(A_q)$ contains the semisimple
part $\xi_s$ of the dilating derivation $\xi:A_q \to A_q$. Since the
endomorphism $\xi_s$ is semisimple and non-trivial, the group $T_q$
is a non-trivial torus. All eigenvectors of the derivation $\xi_s$
in any repreentation of group $GL(A_q)$ are also eigenvectors of the
torus $T_q$; in particular, this applies to the multiplication
element $m \in A_q \otimes A_q \otimes A_q^*$ and to the Poisson
bracket element $p\in A_q \otimes A_q \otimes A_q^*$. Since $\xi_s$
is dilating with non-trivial constant, the torus $T_q$ acts on the
line $k \cdot p \subset A_q \otimes A_q \otimes A_q^*$ by a
non-trivial character $\chi:T \to \gm$.

\begin{lemma}
For every $q > r \geq 1$, the torus $T_q$ preserves the subspace
$\m^r \subset A_q$, and the corresponding reduction map $\red:T_q
\to GL(A_r)$ induces a group isomorphism $T_q \cong T_r$.
\end{lemma}

\proof{} Let $\xi_s,\xi_n \subset \End(A_q)$ be the semisimple and
the nilpotent part of the endomorphism $\xi:A_q \to A_q$. By
\thebull{} the vector field $\xi$ preserves the maximal ideal $\m
\subset A_q$. Since $\codim_k \m = 1$, this is equivalent to
preserving the corresponding line in the dual space
$A_q^*$. Therefore by universality of the Jordan decomposition the
endomorphisms $\xi_s$ and $\xi_n$ also preserve $\m \subset A_k$.
Since by Lemma~\ref{dil} both $\xi_s$ and $\xi_n$ are derivations,
they also preserve the ideal $\m^r \subset A_q$ and act naturally on
the quotient $A_r = A_q/\m^r$.

Denote their reductions $\mod \m^r$ by $\bar\xi_s,\bar\xi_n \in
\End(A_r)$. Then $\bar\xi_s$ is obviously semisimple, $\bar\xi_n$ is
nilpotent, they commute, and $\bar\xi_s + \bar\xi_n$ is the given
derivation $\xi:A_r \to A_r$. By the unicity of the Jordan
decomposition, this means that $\bar\xi_s$ is actually the
semisimple part of the endomorphism $\xi:A_r \to A_r$.

Denote by $P_{qr} \subset GL(A_q)$ the subgroup of endomorphisms
which preserve the ideal $\m^r \subset A_q$, so that we have a
natural reduction map $\red:P_{qr} \to GL(A_r)$, and let $T_{qr} =
\red^{-1}(T_r) \subset P_{qr})$ be the preimage of the torus $T_r
\subset GL(A_r)$ under the natural map $P_{qr} \to GL(A_r)$. Since
$\xi_s:A_q \to A_q$ reduces to $\xi_s:A_r \to A_r$, the Lie algebra
of the subgroup $T_{qr} \subset GL(A_q)$ contains $\xi_s \in
\End(A_q)$. By definition of the group $T_q$ this means $T_q \subset
T_{q,r} \subset P_{qr}$. Therefore $T_q$ indeed preserves $\m^r
\subset A_q$, and the natural reduction map $\red:T_q \to GL(A_r)$
maps $T_q$ into the subgroup $T_r \subset GL(A_r)$.

But the Lie algebra of the image $\red(T_q) \subset T_r$ contains
the endomorphism $\xi_s:A_r \to A_r$. Therefore by definition
$\red(T_q) = T_r$, in other words, the map $\red:T_q \to T_r$ is
surjective.

To prove that this map is injective, it suffices to prove that the
corresponding Lie algebra map is injective. Let $a \in End(A_q)$ be
an element in the Lie algebra of the torus $T_q$ such that $\red(a)
= 0$. Since $T_q$ is a torus, $a$ must be semisimple. On the other
hand, $a$ must be a derivation of the algebra $A_q$ which is zero
$\mod \m^r$ -- in other words, it must send the whole $A_q$ into
$\m^r \subset A_q$. This implies that the endomorphism $a:A_q \to
A_q$ is nilpotent. We conclude that $a=0$.
\endproof

To finish the proof of Proposition~\ref{gr}, it suffices to pass to
the inverse limit. We see that there exists a non-trivial torus $T =
T_q$, $q \geq 1$ which acts on
$$
A = \lim_{\leftarrow}A_q,
$$ 
and a non-trivial character $\chi:T \to \gm = k^*$, such that
$$
t(ab) = t(a)t(b) \qquad
t(\{a,b\}) = \chi(t)\{t(a),t(b)\}
$$
for every $t \in T$, $a,b \in A_\infty$. Take an embedding $\tau:\gm
\to T$ such that $\chi \circ \tau(a) = a^l$ for some non-trivial
integer $l$, and define the grading by means of the induced
$\gm$-action on $A$.
\endproof

\bigskip

\noindent
{\sc Steklov Math Institute\\
Moscow, USSR}

\bigskip

\noindent
{\em E-mail address\/}: {\tt kaledin@mccme.ru}

\end{document}